\documentclass[12pt]{article}
\usepackage{amssymb}
\usepackage{amsmath}
\usepackage{latexsym}
\usepackage{mathrsfs}

\numberwithin{equation}{section}
\newtheorem{thm}[equation]{Theorem}

\newtheorem{lem}[equation]{Lemma}

\title{On a classical spectral optimization problem in linear elasticity\footnote{To appear in the proceedings of the workshop `New Trends in Shape Optimization', Friedrich-Alexander Universit\"{a}t Erlangen-Nuremberg, 23-27 September 2013.}}

\author{Davide Buoso\footnote{corresponding author: dbuoso@math.unipd.it}\ and Pier Domenico Lamberti}

\date{\ }

\begin{document}

\newcommand{\rea}{\mathbb{R}}

\maketitle

%
%
%

\noindent
{\bf Abstract:} We consider a classical shape optimization problem for the eigenvalues of elliptic operators with homogeneous boundary conditions
on domains in the $N$-dimensional Euclidean space. We survey recent results concerning the analytic dependence of the elementary symmetric functions of the eigenvalues
upon domain perturbation and the role of  balls as critical points of such functions subject to volume constraint. Our discussion concerns Dirichlet and buckling-type problems for polyharmonic operators,
the Neumann and the intermediate problems for the biharmonic operator, the Lam\'{e} and the Reissner-Mindlin systems.

\vspace{11pt}


\vspace{6pt}

\section{Introduction}

Let $\Omega$ be a bounded domain (i.e., a bounded connected open set) in ${\mathbb{R}}^N$. As is well known  the problem 
$$
\left\{
\begin{array}{ll}
-\Delta u=\gamma u,& {\rm in}\ \Omega,\\
u=0,& {\rm on}\ \Omega ,  
\end{array}
\right.
$$
admits a divergent sequence of non-negative eigenvalues 
$$
0< \gamma_1[\Omega]< \gamma_2[\Omega]\le \dots \le \gamma_j[\Omega ]\le \dots ,
$$
 where each eigenvalue is repeated as many times as its multiplicity (which is finite).
A classical problem in shape optimization consists in minimizing the eigenvalues $\gamma_j[\Omega ]$ under the assumption that the measure of $\Omega$ is fixed. With regard to this, the most famous 
 result is probably the Rayleigh-Faber-Krahn inequality which reads
\begin{equation}\label{rayleigh}
\gamma_1[\Omega^*]\le \gamma_1 [\Omega],
\end{equation}
where $\Omega^*$ is a ball with the same measure of $\Omega$. In other words, the ball minimizes the first eigenvalue of the Dirichlet Laplacian among all domains with prescribed measure. Note that the first eigenvalue has multiplicity one.    This inequality has been generalized in several directions aiming at  minimization  or maximization results  in the case of  other boundary conditions (for example, Neumann, Robin, Steklov boundary conditions),  other operators  (for example, the biharmonic operator), more general eigenvalue-type problems (for example, the buckling problem for the biharmonic operator)  and other eigenvalues $\gamma_j[\Omega]$ with $j\ne 1$. It is impossible to quote here all available results in this field and we refer to the monographs by Bucur and Buttazzo~\cite{bucbut}, Henrot~\cite{henrot} and Kesavan~\cite{kes} for extensive discussions and references. 

We note that very little is known in the case of polyharmonic operators and systems. We mention that in the case of the biharmonic operator with Dirichlet boundary conditions inequality (\ref{rayleigh}) is known as The Rayleigh Conjecture and has been proved by Nadirashvili~\cite{nad} for $N = 2$ and by Ashbaugh and Benguria~\cite{ash} for $N =2,3$.   We also quote the papers by Bucur, Ferrero and Gazzola~\cite{bucfergaz,bucgaz} concerning the biharmonic operator with Steklov boundary conditions  and Chasman~\cite{chas} for Neumann boundary conditions.  See also the extensive monograph by  Gazzola, Grunau and Sweers~\cite{gaz} for more information on polyharmonic operators.  As for systems, we quote the papers by Kawohl and Sweers~\cite{kaw, kawswe}  which contain interesting lower bounds for the first eigenvalue of the Lam\'{e} system.   

It should be noted that  understanding the behavior of higher eigenvalues is a difficult task even in the case of the Dirichlet Laplacian. A famous result by Buttazzo and Dalmaso~\cite{butdal}  and its recent improvement  by  Mazzoleni and Pratelli~\cite{mazpra} guarantee the existence of a minimizer for $\gamma_j[\Omega]$ in the class of quasiopen sets but no information on the shape of  such minimizer is given. However, it is proved in Wolf and Keller~\cite{wolf} that the minimizers of   higher eigenvalues  in general are not balls and
not even unions of balls. Moreover,  the numerical approach by Oudet~\cite{ou} allows to get an idea of the shape of the minimizers of lower eigenvalues which confirms  the negative result in \cite{wolf}. 

One of the problems arising in the study of higher eigenvalues is related to  bifurcation phenomena associated with the variation of their multiplicity which leads 
to complications such as, for example, lack of differentiability of the eigenvalues with respect to domain perturbation.  However, as it was pointed out in \cite{lala,lalaneu} this problem does not 
affect the elementary symmetric functions of the eigenvalues which depend real-analytically on the domain. This suggests that the elementary symmetric functions of the eigenvalues
might be natural objects in the optimization of multiple eigenvalues. In fact, it turns out that balls are critical points with volume constraint for the elementary symmetric functions of the eigenvalues. 
This property was proved  for  the Dirichlet and Neumann Laplacians  in \cite{lalcri} and later was  proved for polyharmonic operators in \cite{buolam, buolamint}. 

In this survey paper,  we adopt this point of view and show  that the analysis initiated in \cite{lalamult, lala, lalcri, lalaneu} can be extended to a large variety of problems arising in linear elasticity including Dirichlet and buckling-type eigenvalue problems for polyharmonic operators,  biharmonic operator with Neumann and intermediate boundary conditions, Lam\'{e} and Reissner-Mindlin systems.  Details and proofs can be found in  \cite{buoso, buothe, buolam, buolamint, buolamrei,bupro}. 

Our  aim is not only to  collect results spread in different papers but also to present them in a unitary way. In particular,  we provide  a Hadamard-type formula for the shape derivatives of the eigenvalues of the biharmonic operator which is valid not only for Dirichlet boundary conditions (as in the classical case) but also for Neumann and  intermediate boundary conditions. In the case of simple eigenvalues such formula reads
\begin{multline}
\label{hada1}
{\frac{d  \gamma_{n}[\phi_{\epsilon}(\Omega) ]}{d\epsilon}}_{|_{\epsilon =0}}=
\int_{\partial\Omega}\left(|D^2u|^2-2\left(\frac{\partial^2 u}{\partial \nu^2}\right)^2  \right.\\
	\left.+2 \frac{\partial u}{\partial \nu}\left( {\rm div }_{\partial \Omega}[(D^2u)\nu]+\frac{\partial\Delta u}{\partial\nu}\right)
-\gamma u^2 \right)   \psi\cdot n d\sigma,
\end{multline}
where it is assumed that  $\Omega$ is sufficiently smooth, $u$ is an eigenfunction normalized in $L^2(\Omega)$ associated with a simple eigenvalue $\gamma_{n}[ \Omega ]$, and $\phi_{\epsilon}$ are perturbations of the identity $I$ of the type $\phi_{\epsilon}=I+\epsilon \psi$, $\epsilon\in {\mathbb{R}} $. See Theorems~\ref{mainthm}, \ref{nagy} and Lemma~\ref{teclem} for the precise statements and for the case of multiple eigenvalues. Note that in the case of Dirichlet boundary conditions the previous formula gives exactly the 
celebrated Hadamard formula
\begin{equation}
\label{hada2}
{\frac{d  \gamma_{n}[\phi_{\epsilon}(\Omega) ]}{d\epsilon}}_{|_{\epsilon =0}}=-
\int_{\partial\Omega}\left(\frac{\partial^2 u}{\partial \nu^2}\right)^2    \psi\cdot n d\sigma,
\end{equation}
discussed by Hadamard~\cite{had} in the study of a clamped plate  (see also  Grinfeld~\cite{grin}).

This paper is organized as follows: in Section~2 we formulate the eigenvalue problems under consideration, in Section~3 we  state the available analyticity results for  the dependence of the eigenvalues 
upon domain perturbation, in Section~4 we  show that balls are critical points for the elementary symmetric functions of the eigenvalues.

\section{The eigenvalue problems}

Let $\Omega$ be an open set in ${\mathbb{R}}^N$. We denote by $H^m(\Omega)$ the Sobolev space of real-valued functions in $L^2(\Omega)$ with weak derivatives 
up to order $m$ in $L^2(\Omega)$ endowed with its standard norm, and  by $H^m_0(\Omega)$ the closure in $H^m(\Omega)$ of $C^{\infty }_c(\Omega)$.
We consider the following eigenvalue problems on sufficiently regular open sets $\Omega$.   \\

{\bf Dirichlet and buckling problems for polyharmonic operators}

For $m,n\in {\mathbb{N}}$,  $0\le m<n$, we consider the problem 

\begin{equation}
\label{classicnm}
{\mathcal{P}}_{nm}:\ \left\{
\begin{array}{ll}
(-\Delta)^n u=\gamma (-\Delta )^mu, & {\rm in}\ \Omega ,\vspace{2mm}\\
u=\frac{\partial u}{\partial \nu}=\dots =\frac{\partial^{n-1} u}{\partial \nu^{n-1}}=0, & {\rm on }\ \partial \Omega ,
\end{array}\right.
\end{equation}
where  $\nu $ denotes the unit outer normal to $\partial \Omega$.
The case $m=0$ gives the  classical eigenvalue problem for the polyharmonic operator $(-\Delta)^n$ with Dirichlet boundary conditions, while the case $m>0$ represents 
a buckling-type problem. For $N=2$, ${\mathcal{P}}_{10}$ arises for example in the study of a vibrating membrane stretched in a fixed frame, ${\mathcal{P}}_{20}$
corresponds to the case of a vibrating clamped plate and ${\mathcal{P}}_{21}$ is related to plate buckling.   If $\Omega$ is a bounded open set of class $C^1$ then problem 
(\ref{classicnm}) has a sequence of eigenvalues $\gamma^{{\mathcal{P}}_{nm}}_j$ which can be described by the Min-Max Principle. Namely,
\begin{equation}
\label{minmax}
\gamma_j^{{\mathcal{P}}_{nm}}=\min_{\substack{E\subset H^n_0(\Omega ) \\     {\rm dim } E=j}} \max_{\substack{ u\in E \\  u\ne 0}} R_{nm}[u],
\end{equation}
for all $j\in {\mathbb{N}}$, where $R_{nm}[u]$ is the Rayleigh quotient defined by
$$
R_{nm}[u]=\left\{\begin{array}{lll} \frac{  \int_{\Omega }| \Delta^ru|^2 dx  }{ \int_{\Omega }|\Delta^su|^2 dx    },& \ \ {\rm if}\ n=2r,&\ m=2s,\vspace{2.5mm} \\
\frac{  \int_{\Omega }|\Delta^ru|^2 dx  }{ \int_{\Omega }|\nabla \Delta^su|^2 dx    },& \ \ {\rm if}\ n=2r,&\ m=2s+1,\vspace{2.5mm} \\
\frac{  \int_{\Omega }|\nabla \Delta^ru|^2 dx  }{ \int_{\Omega }|\Delta^su|^2 dx    },& \ \ {\rm if}\ n=2r+1,&\ m=2s,\vspace{2.5mm} \\
\frac{  \int_{\Omega }|\nabla \Delta^ru|^2 dx  }{ \int_{\Omega }|\nabla\Delta^su|^2 dx    },& \ \ {\rm if}\ n=2r+1,&\ m=2s+1.    \end{array}  \right.
$$

{\bf  Neumann and intermediate eigenvalue problems for the biharmonic operator}

By Neumann eigenvalue problem for the biharmonic  operator  we mean  the problem

\begin{equation}\label{classicneu}
{\mathcal{N}}:\ \left\{\begin{array}{ll}
\Delta ^2u=\gamma u,\ \ & {\rm in }\ \Omega,\vspace{1mm}\\
\frac{\partial^2u}{\partial^2\nu}=0,\ \ & {\rm on }\ \partial \Omega ,\vspace{1mm} \\
{\rm div}_{\partial\Omega }[ (D^2u) \nu  ]+ \frac{\partial\Delta u}{\partial \nu }=0,\ \ & {\rm on }\ \partial \Omega ,
\end{array}\right.
\end{equation}
where $D^2u$ denotes the Hessian matrix of $u$,  ${\rm div}_{\partial\Omega }$ denotes the tangential divergence operator on  $\partial \Omega$. We recall that ${\rm div}_{\partial\Omega }\, f={\rm div}\, f - [(\nabla f) \nu]\cdot \nu$, for any vector field $f$ smooth enough defined in a neighborhood of $\partial\Omega$.   
Note that we need $\Omega$ to be at least of class $C^2$ for the classical formulation
to make sense, since we need the normal $\nu$ to be differentiable, as can easily
be seen from the boundary conditions; however, we shall interpret problem
(\ref{classicneu}) in the following weak sense 
\begin{equation}\label{weekneu}
\int_{\Omega}D^2u: D^2\varphi dx=\gamma\int_{\Omega}u\varphi dx,\ \ \forall \ \varphi \in  H^2(\Omega),
\end{equation}
where $D^2u: D^2\varphi =\sum_{i,j=1}^Nu_{x_ix_j}\varphi_{x_ix_j}$.
 It is well-known that if $\Omega $ is a bounded open set of class $C^1$ then 
problem 
(\ref{classicneu}) has a sequence of eigenvalues $\gamma^{{\mathcal{N}}}_j$  given by
\begin{equation}
\label{minmaxneu}
\gamma_j^{{\mathcal{N}}}=\min_{\substack{E\subset  H^2(\Omega ) \\     {\rm dim } E=j}} \max_{\substack{ u\in E \\  u\ne 0}} \frac{\int_{\Omega}|D^2u|^2dx}{\int_{\Omega}u^2dx},
\end{equation}
for all $j\in {\mathbb{N}}$, where $|D^2u|^2=\sum_{i,j=1}^Nu^2_{x_ix_j}$. 

If  in (\ref{weekneu}) the space $ H^2(\Omega)$ is replaced by the space $H^2(\Omega)\cap H^1_0(\Omega )$ we obtain the weak formulation of the classical eigenvalue problem
\begin{equation}\label{bihaint}
{\mathcal{I}}:\ \left\{
\begin{array}{ll}
\Delta^2u=\gamma u,\ &\ {\rm in}\ \Omega,\\
u=0,\ &\ {\rm on}\ \partial\Omega,\\
\Delta u -K\frac{\partial u}{\partial \nu}=0,\ &\ {\rm on}\ \partial\Omega,\\
\end{array}
\right.
\end{equation}
where  $K$ denotes the mean curvature of $\partial \Omega $ (the sum of the principal curvatures).
Since   $H^2_{0}(\Omega)\subset H^2(\Omega)\cap H^{1}_0(\Omega )\subset H^{2}(\Omega )$ and the spaces $H^{2}_{0}(\Omega)$, $H^{2}(\Omega)$ are the natural spaces associated with the Dirichlet problem ${\mathcal{P}}_{20}$ and the Neumann problem ${\mathcal{N}}$ respectively, we refer to (\ref{bihaint}) as the 
eigenvalue problem for the biharmonic operator with intermediate boundary conditions. If $\Omega $ is of class $C^1$ 
then problem 
(\ref{bihaint}) has a sequence of eigenvalues $\gamma^{{\mathcal{I}}}_j$  given by
\begin{equation}
\label{minmaxint}
\gamma_j^{{\mathcal{I}}}=\min_{\substack{E\subset H^{2}(\Omega )\cap H^{1}_0(\Omega) \\     {\rm dim } E=j}} \max_{\substack{ u\in E \\  u\ne 0}} \frac{\int_{\Omega}|D^2u|^2dx}{\int_{\Omega}u^2dx},
\end{equation}
for all $j\in {\mathbb{N}}$.\\

{\bf Eigenvalue problem for the Lam\'{e} system}

The eigenvalue problem for the Lam\'{e} system reads

\begin{equation}\label{classiclame}
{\mathcal{L}}:\ \left\{\begin{array}{ll}
-\mu\Delta u -(\lambda+\mu)\nabla {\rm div }u=\gamma   u,\ \ & {\rm in }\ \Omega,\vspace{1mm}\\
u=0,\ \ & {\rm on }\ \partial \Omega ,\vspace{1mm} \\
\end{array}\right.
\end{equation}
where the unknown $u$ is a function taking values in ${\mathbb{R}}^N$ and $\lambda , \mu >0$ are (the Lam\'{e}) constants .  If $\Omega $ is of class $C^1$ then 
 problem 
(\ref{classiclame}) has a sequence of eigenvalues $\gamma^{{\mathcal{L}}}_j$  given by
\begin{equation}
\label{minmaxlame}
\gamma_j^{{\mathcal{L}}}=\min_{\substack{E\subset (H^{1}_0(\Omega ))^N \\     {\rm dim } E=j}} \max_{\substack{ u\in E \\  u\ne 0}} \frac{\int_{\Omega}\mu|\nabla u|^2+(\lambda+\mu){\rm div}^2 udx}{\int_{\Omega}u^2dx},
\end{equation}
for all $j\in {\mathbb{N}}$.\\

{\bf  Eigenvalue problem for the Reissner-Mindlin system}

Finally, we shall consider the eigenvalue problem for the Reissner-Mindlin system 

\begin{equation}\label{reiscla}
	{\mathcal{R}}:\ \left\{
	\begin{array}{l l}
-\frac{\mu}{12}\Delta\beta-\frac{\mu+\lambda}{12}\nabla\mathrm{div}\beta  -\mu\frac{\kappa}{t^2}(\nabla w-\beta)
		=\frac{t^2\gamma}{12}\beta, & \mathrm{in}\ \Omega	,\vspace{0,2cm}\\
	-\mu\frac{k}{t^2}(\Delta w -\mathrm{div}\beta)=\gamma w, & \mathrm{in}\ \Omega
 , \vspace{0,2cm}\\
	\beta =0, \ w=0, & \mathrm{on}\ \Omega ,
	\end{array}
	\right.
	\end{equation}
where the unknown $(\beta , w)=(\beta_1, \dots , \beta_N, w)$ is a function with values in $\mathbb{R}^{N+1}$ and $\lambda, \mu , \kappa, t> 0$ are constants.  According to the Reissner-Mindlin model for moderately thin plates, for $N=2$ system (\ref{reiscla}) describes the free vibration modes of an elastic clamped
plate $\Omega \times (-t/2,t/2)$ with midplane $\Omega$ and thickness $t$. In that  case  
  $\lambda$ and $\mu$ are the Lam\'{e} constants, $\kappa$ is the correction factor, $w$ the transverse displacement of the midplane, $\beta =(\beta_1,\beta_2)$   
the corresponding rotation and $t^2\gamma$ the vibration frequency.

If $\Omega $ is of class $C^1$ 
then problem 
(\ref{reiscla}) has a sequence of eigenvalues $\gamma^{{\mathcal{R}}}_j$  given by
\begin{equation}
\label{minmaxreis}
\gamma_j^{{\mathcal{R}}}=\min_{\substack{E\subset (H^{1}_0(\Omega ))^{N+1} \\   {\rm dim } E=j}} \max_{\substack{ (\beta , w)\in E \\  u\ne 0}} \frac{ \int_{\Omega}\frac{\mu}{12} |\nabla\beta |^2\eta +\frac{\mu+\lambda}{12} \mathrm{div}^2\beta  +\mu\frac{\kappa}{t^2} |\nabla w -\beta |^2dx        }{ \int_{\Omega}w^2+\frac{t^2}{12}|\beta |^2 dx  },
\end{equation}
for all $j\in {\mathbb{N}}$.\\

\section{Analyticity results }

Let $\Omega $ be a bounded open set in ${\mathbb{R}}^N$ of class $C^1$.  In the sequel, we shall consider problems (\ref{classicnm}),   (\ref{classicneu}), (\ref{bihaint}), (\ref{classiclame}), (\ref{reiscla})   on  families of open sets  parametrized by suitable diffeomorphisms $\phi $
defined on $\Omega $. To do so, for $k\in {\mathbb{N}}$ we set
$$
{\mathcal{A}}_{\Omega }^k=\biggl\{\phi\in C^k_b(\Omega\, ; {\mathbb{R}}^N ):\ \inf_{\substack{x_1,x_2\in \Omega \\ x_1\ne x_2}}\frac{|\phi(x_1)-\phi(x_2)|}{|x_1-x_2|}>0 \biggr\},
$$
where $C^k_b(\Omega\, ; {\mathbb{R}}^N )$ denotes the space of all functions from $\Omega $ to ${\mathbb{R}}^N$ of class $C^k$, with bounded derivatives up to order
  $k$.  Note that if $\phi \in {\mathcal{A}}_{\Omega }^k$ then $\phi $ is injective, Lipschitz continuous and $\inf_{\Omega }|{\rm det }\nabla \phi |>0$. Moreover, $\phi (\Omega )$ is a bounded open set of class $C^1$ and the inverse map $\phi^{(-1)}$ belongs to  ${\mathcal{A}}_{\phi(\Omega )}^k$.   
Thus it is natural to consider the above mentioned eigenvalue problems  on $\phi (\Omega )$ and study  the dependence of the corresponding eigenvalues $\gamma_j^{\mathcal{P}_{nm}}[\phi (\Omega )]$, $\gamma_j^{\mathcal{N}}[\phi (\Omega )]$, $\gamma_j^{\mathcal{I}}[\phi (\Omega )]$, $\gamma_j^{\mathcal{L}}[\phi (\Omega )]$, 
$\gamma_j^{\mathcal{R}}[\phi (\Omega )]$ on $\phi \in {\mathcal{A}}_{\Omega }^k$ for suitable values of $k$. \\
The choice of $k$ depends on the problem.  In the sequel it will be always understood that $k$ is chosen as follows:\\
\begin{equation}\label{kappach}
\begin{array}{ll}
{\rm  Problem}\ {\mathcal{P}}_{nm}:& k=n,\\
{\rm  Probelms}\ {\mathcal{N}}\ {\rm and}\ {\mathcal{I}}:& k=2,\\
{\rm  Problems}\  {\mathcal{L}}\  {\rm and}\ {\mathcal{R}}:&  k=1.\\
\end{array}
\end{equation}
Moreover,  we shall simply write  $\gamma_j[\phi ]$ instead of $\gamma_j^{\mathcal{P}_{nm}}[\phi (\Omega )]$, $\gamma_j^{\mathcal{N}}[\phi (\Omega )]$, $\gamma_j^{\mathcal{I}}[\phi (\Omega )]$, $\gamma_j^{\mathcal{L}}[\phi (\Omega )]$, 
$\gamma_j^{\mathcal{R}}[\phi (\Omega )]$, with the understanding that our statements refer to any of the problems (\ref{classicnm}),   (\ref{classicneu}), (\ref{bihaint}), (\ref{classiclame}), (\ref{reiscla}).\\
We endow the space $C^k_b(\Omega\, ; {\mathbb{R}}^N )$ with its usual norm defined by $\| f \|_{C^k_b(\Omega\, ;{\mathbb{R}}^N )}=\sup_{|\alpha|\le k,\ x\in\Omega  } |D^{\alpha }f(x)|$. We recall  that ${\mathcal{A}}_{\Omega }^k$ is an open set in  $C^k_b(\Omega\, ;{\mathbb{R}}^N )$, see \cite[Lemma~3.11]{lala}.
    Thus, it makes sense to  study differentiability and analyticity properties of the maps $\phi \mapsto \gamma_j[\phi (\Omega )]$ defined for $\phi \in {\mathcal{A}}_{\Omega }^k$.
  
As in \cite{lala}, we fix a finite set of indexes $F\subset \mathbb{N}$
and we consider those maps $\phi\in {\mathcal{A}}^k_{\Omega }$ for which the eigenvalues
with indexes in $F$  do not coincide with eigenvalues with indexes not
in $F$; namely we set
$$
{\mathcal { A}}^{k}_{F, \Omega }= \left\{\phi \in {\mathcal { A}}^k_{\Omega }:\
\gamma_j[\phi ]\ne \gamma_l[\phi],\ \forall\  j\in F,\,   l\in \mathbb{N}\setminus F
\right\}.
$$
It is also convenient to consider those maps $\phi \in {\mathcal { A}}^{k}_{F, \Omega } $ such that all the eigenvalues with index in $F$
 coincide and set
$$
\Theta_{F, \Omega }^{k} = \left\{\phi\in {\mathcal { A}}_{F, \Omega }^{k}:\ \gamma_{j_1}[\phi ]
=\gamma_{j_2}[\phi ],\, \
\forall\ j_1,j_2\in F  \right\} .
$$
For $\phi \in {\mathcal { A}}^{k}_{F, \Omega }$, the elementary symmetric functions of the eigenvalues with index in $F$ are defined by
\begin{equation}
\label{sym1}
\Gamma_{F,h}[\phi ]=\sum_{ \substack{ j_1,\dots ,j_h\in F\\ j_1<\dots <j_h} }
\gamma_{j_1}[\phi ]\cdots \gamma_{j_h}[\phi ],\ \ \ h=1,\dots , |F|.
\end{equation}

In order to state Theorems~\ref{mainthm}, \ref{nagy}, we need to define a quantity $M[u,v]$ where  $u,v$ are eigenfunctions associated with an eigenvalue $\gamma$ on a  smooth bounded open set $\Omega$.  For each problem, $M[u,v]$ is a real valued function  defined on $\partial \Omega$ as follows:

\begin{itemize}
\item  Problem ${\mathcal{P}}_{nm}$:
\begin{equation}\label{dirhad}
M[u,v]=\frac{\partial^n u}{\partial \nu^n}\frac{\partial^n v}{\partial \nu^n};
\end{equation}
\item  Problem ${\mathcal{N}}$:
\begin{equation}\label{neuhad}
M[u,v]=\gamma uv -D^2u : D^2v ;
\end{equation}
\item  Problem ${\mathcal{I}}$:
\begin{equation}\label{inthad}
M[u,v]=D^2u: D^2v -  2\Delta_{\partial\Omega}\biggl(\frac{\partial u}{\partial\nu}\frac{\partial v}{\partial\nu}\biggr)
				-\left(\frac{\partial u}{\partial\nu}\frac{\partial^3v}{\partial\nu^3}
					+\frac{\partial u}{\partial\nu}\frac{\partial^3v}{\partial\nu^3}\right)      ;
\end{equation}
\item  Problem ${\mathcal{L}}$:
\begin{equation}
M[u,v]= \mu\frac{\partial u}{\partial \nu}\cdot\frac{\partial v}{\partial \nu}
	+(\mu+\lambda)\left(\frac{\partial u}{\partial \nu}\cdot \nu\right)\left(\frac{\partial v}{\partial \nu}\cdot \nu\right); 
	\end{equation}

\item  Problem ${\mathcal{R}}$:
\begin{equation}
M[u,v]	=  \frac{\mu}{12}  \frac{\partial\beta}{\partial \nu}\cdot\frac{\partial\theta}{\partial \nu}+\frac{\mu+\lambda}{12}\left(\frac{\partial\beta}{\partial \nu}\cdot \nu\right)\left(\frac{\partial\theta}{\partial \nu}\cdot \nu	\right)+\frac{\kappa\mu}{t^2}\frac{\partial w}{\partial \nu}\frac{\partial u}{\partial \nu}  ;
\end{equation}
where $u=(\beta , w)$ and $v=(\theta , u )$. 
\end{itemize}

In (\ref{inthad}), $\Delta_{\partial\Omega}$ denotes the tangential Laplacian on $\partial \Omega$. Recall that $\Delta_{\partial\Omega}u={\rm div}_{\partial \Omega }\nabla_{\partial \Omega }u$ where $\nabla_{\partial \Omega }u=\nabla u- \frac{\partial u}{\partial \nu}\nu$ is the tangential gradient of $u$.

Moreover, formula (\ref{derivn}) below is expressed in terms of  a basis $\{u_l\}_{F} $ of the eigenspace associated with an eigenvalue $\gamma$ on an open set $\tilde \phi (\Omega)$. It will be understood that such basis is orthonormal with respect to the appropriate  $L^2$-scalar product, which is the standard scalar product in $L^2(\tilde\phi (\Omega ))$ for problems (\ref{classicnm}) with $m=0$,   (\ref{classicneu}), (\ref{bihaint}), (\ref{classiclame}) and the scalar product defined by $\int_{\tilde\phi (\Omega)}(wv+\frac{t^2}{12}\beta\cdot\eta)dy  $ for problem (\ref{reiscla}).  Note that in the case of problem (\ref{classicnm}) with arbitrary $m$, we use the natural scalar product  associated with right-hand side of the equation, i.e., the scalar product defined by $\int_{\tilde\phi (\Omega)}\Delta^{\frac{m}{2}}u\Delta^{\frac{m}{2}}vdy$ if $m$ is even, and   $\int_{\tilde\phi (\Omega)}\nabla \Delta^{\frac{m-1}{2}}u\nabla \Delta^{\frac{m-1}{2}}vdy$ if $m$ is odd.

Then we have the following

\begin{thm}\label{mainthm}Let $\Omega $ be a bounded open set in ${\mathbb{R}}^N$ of class $C^1$  and  $F$ be a finite set in  ${\mathbb{N}}$. Let $k\in {\mathbb{N}}$ be as in (\ref{kappach}).  The set ${\mathcal { A}}^{k}_{F, \Omega }$ is open in
$C^k_b(\Omega\, ; {\mathbb{R}}^N )$ and the real-valued maps which take $\phi\in {\mathcal { A}}^{k}_{F, \Omega } $ to $ \Gamma_{F,h}[\phi ]$ are real-analytic on  ${\mathcal { A}}^{k}_{F, \Omega }$ for all $h=1,\dots , |F|$. Moreover, if $\tilde \phi\in \Theta^{k}_{F, \Omega }  $ is such that the eigenvalues $\gamma_j[\tilde \phi]$ assume the common value $\gamma_F[\tilde \phi ]$ for all $j\in F$, and  $\tilde \phi (\Omega )$ is of class $C^{2k}$ then  the Frech\'{e}t differential of the map $\Gamma_{F,h}$ at the point $\tilde\phi $ is delivered by the formula
\begin{equation}
		\label{derivn}
			d|_{\phi=\tilde{\phi}}\Gamma_{F,h}[\psi]=-\gamma_F^{h-1}[\tilde{\phi}]\binom{|F|-1}{h-1}
			\sum_{l=1}^{|F|}\int_{\partial \tilde\phi (\Omega ) } M[u_l,u_l]\zeta \cdot \nu d\sigma  ,
		\end{equation}
		for all $\psi\in C^{k}_b(\Omega ; {\mathbb{R}}^N)$, where $\{u_l\}_{l\in F}$ is an orthonormal basis  of the eigenspace associated with $\gamma_F[\tilde \phi]$, and $\zeta =\psi\circ\tilde{\phi}^{(-1)}$.
\end{thm}
The proof of this theorem can be done by adapting that of \cite[Theorem 3.38]{lala} (see also \cite[Theorem 2.5]{lalaneu}). Namely, by pulling-back to $\Omega$ via $\phi$ the operator defined on $\phi (\Omega)$, one reduces the problem to the study of a family of operators $T_{\phi}$ defined on the fixed domain $\Omega$. Such operators turn out to be self-adjoint with respect to a scalar product also depending on $\phi$, which is obtained by pulling-back the appropriate scalar product defined of $L^2(\phi (\Omega))$. Then it is possible to apply the abstract results of  \cite{lala} in order to prove the real-analyticity of the symmetric functions of the eigenvalues. Formula (\ref{derivn}) is also deduced by a general formula concerning the eigenvalues of self-adjoint operators proved in \cite{lala} combined with lenghty calculations which depend on the specific case under consideration. We refer to the papers indicated in the introduction for details. \\

If we consider 
 domain perturbations depending real analytically on one scalar parameter, it is possible to describe  all the eigenvalues
splitting from a multiple eigenvalue of multiplicity $m$ by means of $m$ real-analytic functions. For the sake of completeness we state the  following Rellich-Nagy-type theorem which can be proved by using the abstract  results \cite[Theorem~2.27, Corollary~2.28]{lala} which, in turn, are proved by an argument based on reduction to finite dimension.

\begin{thm}\label{nagy}Let $\Omega$ be a bounded open set in ${\mathbb{R}}^N$ of class $C^1$.  Let $k\in {\mathbb{N}}$ be as in (\ref{kappach}), $\tilde \phi \in {\mathcal{A}}_{\Omega}^k$ and $\{\phi_{\epsilon}\}_{\epsilon\in {\mathbb{R}}} \subset {\mathcal{A}}_{\Omega}^k$ be a family depending real-analytically on $\epsilon$ such that $\phi_0=\tilde \phi$. Let $\tilde\gamma$ be an eigenvalue on $\tilde\phi (\Omega)$ of multiplicity $m$, namely $\tilde \gamma =\gamma_{n,t}[\tilde \phi ]=\dots =\gamma_{n+m-1,t}[\tilde \phi]$ for some $n\in {\mathbb{N}}$. Then there exists an open interval $I$ containing zero and $m$ real-analytic functions $g_{1},\dots , g_m$ from $I$ to ${\mathbb{R}}$ such that $\{\gamma_{n,t}[ \phi_{\epsilon} ],\dots ,\gamma_{n+m-1,t}[\phi_{\epsilon}]\}=\{g_1(\epsilon ), \dots  , g_m(\epsilon) \}$ for all $\epsilon \in I$. Moreover, if $\tilde \phi (\Omega)$  is an open set of class $C^{2k}$ then the derivatives 
 $g'_1(0), \dots , g_m'(0)$ at zero of the functions $g_1, \dots , g_m$ coincide with the eigenvalues of the matrix 
$$\left(- \int_{\tilde\phi (\Omega)}M[u_i,u_j] \zeta\cdot \nu d\sigma  \right)_{i,j\in \{1, \dots , m\}}$$ where $u_i$, $i=1,\dots ,m$, is an orthonormal basis  of the eigenspace associated with $\tilde\gamma$. 
\end{thm}

In the case of the biharmonic operator the quantities $M[u_i,u_j]$ can be represented by one single formula which is valid for the Dirichlet problem ${\mathcal{P}}_{20}$, the Neumann problem 
${\mathcal{N}}$ and the intermediate problem ${\mathcal{I}}$. Indeed, we can prove the following 

\begin{lem}\label{teclem}Let  $u,v$  be eigenfunctions associated with the same eigenvalue $\gamma$ of one of the problems ${\mathcal{P}}_{20}$, ${\mathcal{N}}$, ${\mathcal{I}}$ on a  bounded open set $\Omega$ of class $C^4$. Then 
\begin{eqnarray}\label{unique}
M[u,v]&=&2\frac{\partial^2 u}{\partial \nu^2}\frac{\partial^2 v}{\partial \nu^2}   -D^2u: D^2v +\gamma uv -\frac{\partial u}{\partial \nu}\left( {\rm div }_{\partial \Omega}[(D^2v)\nu]+ \frac{\partial \Delta v}{\partial \nu }  \right) \nonumber  \\
& &-\frac{\partial v}{\partial \nu}\left( {\rm div }_{\partial \Omega}[(D^2u)\nu]+ \frac{\partial \Delta u}{\partial \nu }  \right).
\end{eqnarray} 
In particular, in these cases formula (\ref{derivn}) reads 
\begin{eqnarray}
		\label{derivn2}
\lefteqn{
			d|_{\phi=\tilde{\phi}}\Gamma_{F,h}[\psi]=-\gamma_F^{h-1}[\tilde{\phi}]\binom{|F|-1}{h-1}
			\sum_{l=1}^{|F|}\int_{\partial \tilde\phi (\Omega ) } \left(  2\left(\frac{\partial^2 u_l}{\partial \nu^2}\right)^2-  |D^2u_l|^2  \right. }\nonumber \\
& & \qquad\qquad\qquad \left. +\gamma u_l^2-2 \frac{\partial u_l}{\partial \nu}\left( {\rm div }_{\partial \Omega}[(D^2u_l)\nu]+\frac{\partial\Delta u_l}{\partial\nu}\right)\right) \zeta \cdot \nu d\sigma  .
\end{eqnarray}
\end{lem}

{\bf Proof.} In the case of problem ${\mathcal{P}}_{20}$, taking into account the boundary conditions $u=v=0$ on $\partial \Omega $ and $ \nabla u=\nabla v =0$ on $\partial \Omega $, it follows that 
$D^2u: D^2v=\frac{\partial^2 u}{\partial \nu^2}\frac{\partial^2 v}{\partial \nu^2}  $ on $\partial \Omega $ hence 
the right-hand side of (\ref{unique}) equals  the right-hand side of (\ref{dirhad}) with $n=2$. 

In the case of problem ${\mathcal{N}}$, functions $u$ and $v$ satisfy the boundary conditions in (\ref{classicneu}) hence we  immediately conclude that the right-hand side of (\ref{unique}) equals  the right-hand side of (\ref{neuhad}). 

Finally, we consider the intermediate problem ${\mathcal{I}}$. In this case, several calculations are required. To begin with, we note that since $u=v=0$ on $\partial \Omega $ we have
\begin{eqnarray}\label{contiint1}\lefteqn{
\Delta_{\partial \Omega }\left(\frac{\partial u}{\partial \nu}\frac{\partial v}{\partial \nu}   \right)= \Delta_{\partial \Omega }\left(\frac{\partial u}{\partial \nu} \right)\frac{\partial v}{\partial \nu}  +2\nabla_{\partial \Omega }\frac{\partial u}{\partial \nu}  \nabla_{\partial \Omega}\frac{\partial v}{\partial \nu}
+ \frac{\partial u}{\partial \nu} \Delta_{\partial \Omega }\left(\frac{\partial v}{\partial \nu} \right)}\nonumber \\
& & \qquad= {\rm div }_{\partial \Omega}[(D^2u)\nu] \frac{\partial v}{\partial \nu}  +2\nabla_{\partial \Omega }\frac{\partial u}{\partial \nu}  \nabla_{\partial \Omega}\frac{\partial v}{\partial \nu}
+ \frac{\partial u}{\partial \nu}{\rm div }_{\partial \Omega}[(D^2v)\nu].
\end{eqnarray}
On the other hand, we have
\begin{eqnarray}\label{contiint2}\lefteqn{
\Delta_{\partial \Omega }\left(\frac{\partial u}{\partial \nu}\frac{\partial v}{\partial \nu}   \right)=
 \Delta_{\partial \Omega }\left(\nabla u\cdot \nabla v   \right)=  \Delta \left(\nabla u\cdot \nabla v   \right)-\frac{\partial^2 (\nabla u\cdot \nabla v) }{\partial \nu^2}
-K \frac{\partial (\nabla u\cdot \nabla v) }{\partial \nu} }\nonumber\\
& &= \nabla \Delta u \cdot \nabla v  +\nabla \Delta v \cdot \nabla u +2D^2u: D^2v -2[(D^2u)\nu)]\cdot [(D^2v)\nu]-\nabla u\frac{\partial^2\nabla v}{\partial \nu^2}\nonumber\\
& & -\nabla v\frac{\partial^2\nabla u}{\partial \nu^2}-K\nabla u \cdot \frac{\partial \nabla v }{\partial \nu }-K\nabla v \cdot \frac{\partial \nabla u }{\partial \nu }=\frac{\partial \Delta u}{\partial \nu}\frac{\partial v}{\partial \nu }  + \frac{\partial \Delta v}{\partial \nu}\frac{\partial u}{\partial \nu }\nonumber \\ 
& & +2D^2u: D^2v    - 2\nabla_{\partial \Omega }\frac{\partial u}{\partial \nu}  \nabla_{\partial \Omega}\frac{\partial v}{\partial \nu}
-\frac{\partial u}{\partial \nu} \frac{\partial^3 v}{\partial \nu^3}
-\frac{\partial v}{\partial \nu} \frac{\partial^3u}{\partial \nu^3}-K \frac{\partial u}{\partial \nu} \frac{\partial^2 v}{\partial \nu^2}\nonumber \\
& & 
-K \frac{\partial v}{\partial \nu} \frac{\partial^2 u}{\partial \nu^2}.
\end{eqnarray}
By taking into account that functions $u$ and $v$ satisfy the boundary condition $\frac{\partial^2 u}{\partial \nu^2}=\frac{\partial^2 v}{\partial \nu^2}=0$ on $\partial \Omega $, and by summing the first and last terms in the respective equalities (\ref{contiint1}) and (\ref{contiint2}) we get
\begin{eqnarray}\lefteqn{
2\Delta_{\partial \Omega }\left(\frac{\partial u}{\partial \nu}\frac{\partial v}{\partial \nu}   \right)=2D^2u: D^2v +  \frac{\partial u}{\partial \nu}\left( {\rm div }_{\partial \Omega}[(D^2v)\nu]+ \frac{\partial \Delta v}{\partial \nu }  \right)}\nonumber \\
& &    \frac{\partial v}{\partial \nu}\left( {\rm div }_{\partial \Omega}[(D^2u)\nu]+ \frac{\partial \Delta u}{\partial \nu }  \right)  -\frac{\partial u}{\partial \nu} \frac{\partial^3 v}{\partial \nu^3}-\frac{\partial v}{\partial \nu} \frac{\partial^3 u}{\partial \nu^3}.
\end{eqnarray}
The previous equality shows that in the case of problem $\mathcal{I}$   the right-hand side of (\ref{unique}) equals  the right-hand side of (\ref{inthad}). \hfill $\Box$

\section{Isovulmetric perturbations}

Consider the following extremum problems for the symmetric functions of the eigenvalues
\begin{equation}\label{min}
\min_{V[\phi ]={\rm const}}\Gamma_{F,h}[\phi ]\ \ \ {\rm or}\ \ \ \max_{   V[\phi ]={\rm const}}\Gamma_{F,h}[\phi],
\end{equation}
where $V[\phi ]$ denotes the $N$-dimensional Lebesgue measure of $\phi (\Omega )$.

Note that if $\tilde \phi \in {\mathcal{A}}_{\Omega }^k$ is a minimizer or maximizer in (\ref{min}) then $\tilde \phi $ is a critical domain transformation for the map $\phi\mapsto  \Gamma_{F,h}[\phi ]$ subject to volume constraint, {\rm i.e.,}
\begin{equation}\label{inc}
{\rm Ker\ d}|_{\phi =\tilde \phi }V \subset {\rm Ker\ d}|_{\phi =\tilde \phi } \Lambda_{F,h},
\end{equation}
where $V$ is the real valued function defined on ${\mathcal{A}}_{\Omega }^k$ which takes $\phi \in {\mathcal{A}}_{\Omega }^k$ to $V[\phi ]$.

The following theorem provides a characterization of all critical domain transformations $\phi$. See \cite{lalcri}
for the case of the Dirichlet and Neumann Laplacians.

\begin{thm}\label{car}Let $\Omega $ be a bounded open set in ${\mathbb{R}}^N$ of class $C^1$.  Let $k\in {\mathbb{N}}$ be as in (\ref{kappach}). Let $F$ be a finite subset of ${\mathbb{N}}$.
Assume that $\tilde \phi\in \Theta_{F, \Omega }^{k} $ is such that $\tilde \phi (\Omega )$ is of class $C^{2k}$ and that the eigenvalues $\gamma_j[\tilde \phi ]$ have the common value $\gamma_F[\tilde \phi ]$
for all $j\in F$. Let $\{  u_l\}_{l\in F}$ be an orthornormal basis  of the eigenspace corresponding to $\gamma_F[\tilde \phi ]$. Then  $\tilde \phi$
is a critical domain transformation  for any of the functions $\Gamma_{F,h}$, $h=1,\dots , |F|$,  with  volume constraint
if and only if  there exists $C \in {\mathbb{R}}$ such that 
\begin{equation}\label{ovcumu}
\sum_{l\in F} M[u_l,u_l]= C,\ \ {\rm on}\ \partial\tilde \phi (\Omega ) . 
\end{equation}
\end{thm}
Formula (\ref{ovcumu}) follows from an application of the Lagrange Multipliers Theorem (see  e.g., Deimling~\cite[\S~26]{dei} for a formulation valid in the case of infinite dimensional spaces) and formula (\ref{derivn}).\\

Finally, thanks to the rotation invariance of the operators related to the problems we have considered, it is possible to prove the following 

\begin{thm}\label{ballthm}
Let the same assumptions of Theorem~\ref{car} hold. If $\tilde \phi (\Omega )$ is a ball then condition  (\ref{ovcumu}) is satisfied.
\end{thm}
The proof of this theorem is based on the following main idea.  First, we assume that $\tilde \phi (\Omega )$ is a ball with radius $R$ centered at zero.  In the case of polyharmonic operators, we have that by the rotation invariance of the Laplace operator, if $\{u_l\}_{l\in F}$ is an orthonormal basis of an eigenspace, then  $\{ u_{l}\circ A\}_{l\in F}$   is also an orthonormal basis of the same  eigenspace  for all $A\in O_{N}({\mathbb{R}})$, where  $O_{N}({\mathbb{R}})$ denotes the group of orthogonal linear transformations in ${\mathbb{R}}^N$. 
 Since both $\{ u_l\}_{l\in F}$
and $\{  u_{l}\circ A\}_{l\in F}$ are orthonormal bases of the same space, it follows that
$
\sum_{l=1}^{|F|} { u}_{l}^{2}\circ
A=\sum_{l=1}^{|F|} { u}_{l}^{2}\, ,
$
 for all $A\in O_{N}({\mathbb{R}})$. Thus $\sum_{l=1}^{|F|}{ u}_{l}^{2}$ is a radial function.  Then the radiality of $\sum_{l=1}^{|F|}{ u}_{l}^{2}$ combined with appropriate calculations and similar arguments as above, allows to conclude that 
(\ref{ovcumu}) is satisfied. (Note that in the case of vector-valued functions, say in the case of the Lam\'{e} system for simplicity, one has clearly to rotate the vector itself by considering $A^{t}\cdot (u_{l}\circ A)$ where we identify $A$ with its matrix.) \\

It would be interesting to describe the family of open sets   $\tilde \phi (\Omega)$  for which condition  (\ref{ovcumu}) is satisfied.   In the case of problem ${\mathcal{P}}_{10}$ a classical result by Serrin~\cite{serrin} guarantees that  if condition  (\ref{ovcumu}) is satisfied for the first eigenfunction then  $\tilde \phi (\Omega )$ is a ball. The same result has been proved by Dalmasso~\cite{dalmasso} in the case of problem ${\mathcal{P}}_{20}$
 under the assumption that the first eigenfunction does not change sign; for problem ${\mathcal{P}}_{21}$ a different method by Weinberger and Willms leads to the same conclusion (see e.g., \cite{henrot}).\\

{\bf Acknowledgments}: The  authors acknowledge financial support from the research project  `Singular perturbation problems for differential operators',  Progetto di Ateneo of the University of Padova. The authors are members of the Gruppo Nazionale per
l'Analisi Matematica, la Probabilit\`{a} e le loro Applicazioni (GNAMPA)
of the Istituto Nazionale di Alta Matematica (INdAM).

\noindent {\small
Davide Buoso$^{1}$ and Pier Domenico Lamberti$^{2}$\\
Dipartimento di Matematica\\
Universit\`{a} degli Studi di Padova\\
Via Trieste,  63\\
35126 Padova\\
Italy\\
e-mail$^{1}$:	dbuoso@math.unipd.it \\
e-mail$^{2}$:	lamberti@math.unipd.it}

\end{document}